\documentclass[12pt]{amsart}

\usepackage{amscd,amssymb,amsfonts,amsmath,amstext,amsthm}

\newcommand{\ot}{\otimes}

\newcommand{\eps}{\varepsilon}
\newcommand{\al}{\alpha}
\newcommand{\La}{\Lambda}
\newcommand{\vs}{\vskip 3pt}
\newcommand{\ds}{\displaystyle}

\def\CO{{\mathcal O}}

\def\BZ{{\mathbb Z}}
\def\BC{{\mathbb C}}
\def\BR{{\mathbb R}}

\newtheorem{thm}{Theorem}[section]
\newtheorem{df}[thm]{Definition}
\newtheorem{remark}[thm]{Remark}
\newtheorem{ex}[thm]{Example}
\newtheorem{cor}[thm]{Corollary}
\newtheorem{prop}[thm]{Proposition}
\newtheorem{lem}[thm]{Lemma}

\newcommand{\num}{\refstepcounter{thm}}

\title[Representations of Hopf Algebras associated to $S_n$]
{Representations of some Hopf Algebras associated to the symmetric group $S_n$}
\author{Andrea Jedwab}
\address{University of Southern California, Los Angeles, CA 90089-1113}
\email{jedwab@usc.edu}
\author{Susan Montgomery}
\address{University of Southern California, Los Angeles, CA 90089-1113}
\email{smontgom@math.usc.edu}

\begin{document}

\maketitle


\section{Introduction}

In this paper we study the representations of two semisimple Hopf algebras related to the symmetric 
group $S_n$, namely the bismash products $H_n = k^{C_n}\# kS_{n-1}$ and its dual 
$J_n = k^{S_{n-1}}\# kC_n = (H_n)^*,$ where $k$ is an algebraically closed field of characteristic 0. 
Both algebras are constructed using the standard representation of $S_n$ as a factorizable group, that is 
$S_n = S_{n-1}C_n = C_nS_{n-1}$, where (the copy of ) $C_n$ is generated by $(1,2,3,\ldots, n)$. We are 
interested especially in determining the Frobenius-Schur indicators of the simple modules for 
$H_n$ and $J_n$. We prove that for $H_n$, the indicators of all simple modules are +1, that is, the 
algebra is {\it totally orthogonal}. This fact was known classically for $S_n$ itself, and was also shown 
to be true for the Drinfeld double $D(S_n)$ of  $S_n$ in \cite{KMM}.   

However, for the dual Hopf algebra $J_n = k^{S_{n_1}}\# kC_n$, the situation is considerably more 
complicated: the indicator can have value 0 as well as 1. When $n = p$, a prime, we obtain a precise result 
as to which representations have indicator +1 and which ones have 0; in fact as $p \to \infty$, the proportion of 
simple modules with indicator 1 becomes arbitrarily small. To prove this, we first prove a result about 
Frobenius-Schur indicators for more general bismash products $H =k^G\# kF$, coming from any factorizable 
group of the form $L = FG$ such that $F\cong C_p.$

It was shown in \cite{LM} that the classical theorem of Frobenius and Schur, giving a 
formula to compute the indicator of a simple module $V$ for a finite group $G,$ extends 
to any semisimple Hopf algebra $H$. This fact was shown a bit earlier for the special case of Kac 
algebras over $\BC$ in \cite{FGSV}. For such an $H$-module $V$, with character $\chi$ and indicator 
$\nu(V) = \nu (\chi)$, the only possible values of $\nu (\chi)$ are $0,1$ and $-1$. 
Then $\nu(V)\ne 0$ if and only if $V$ is self dual; assuming $V$ is self-dual, then $\nu(V)=+1$ 
if and only if $V$ admits a non-degenerate $H$-invariant symmetric bilinear form and 
$\nu(V)=-1$ if and only if $V$ admits a non-degenerate $H$-invariant skew-symmetric bilinear
form. 

However it is not true for Hopf algebras over $\BC$, as it was for the case of a finite group $G$, that 
$\nu(V)=+1$ if and only if $V$ is defined over $\BR$. 

The case when $\nu(V) = +1$ for all simple $G$-modules $V$ has been of particular interest; 
such groups are called {\it totally orthogonal} in \cite{GW}. This terminology seems suitable for Hopf 
algebras as well, in view of the existence of the bilinear forms described above. In particular it was
known classically that if $G$ is any finite real reflection group, then all $\nu(V) = +1$ \cite{H}; this 
includes the case $G = S_n$, as noted above. Some other examples are given in \cite{GW}. 

Recently it was shown in \cite{GM} that $D(G)$, the Drinfel'd double of a finite group, is totally orthogonal 
for any finite real reflection group $G$. 

We remark that the indicator is proving to be a very useful invariant in Hopf algebras. It has been used 
in the classification of Hopf algebras, in finding the dimensions of simple modules \cite{KSZ1}, and in 
determining the prime divisors of the exponent of a (semisimple) Hopf algebra \cite{KSZ2}. Moreover, 
the indicator is an invariant of monoidal categories \cite{MaN}. Finally, bismash products constructed from 
factorizable groups, and bicrossed products, a generalization of bismash products using group cocycles, 
are among the most fundamental constructions of Hopf algebras. As of this writing, there are no known 
semisimple Hopf algebras over $\BC$ other than (Drinfel'd) twists of the quasi-Hopf version of bicrossed 
products \cite{N1}. 

The particular bismash products $H_n = k^{C_n}\# kS_{n-1}$ and $J_n = k^{S_{n-1}}\# kC_n$
which we study here have several additional nice properties. First, any bicrossed product
constructed from the factorization $S_n = S_{n-1}C_n = C_nS_{n-1}$ gives nothing new;
that is, it is isomorphic to $H_n$ or $J_n$ \cite{Ma2}; see Example \ref{Sn}. Second,
it has recently been shown by Collins that the algebra structure of $H_n$ is unlike that
of any group algebra; that is, there does not exist any finite group $G$ such that the set of degrees of
the irreducible representations of $G$ is the same as the set of degrees of the irreducible
representations of $H_n$, for any prime $p>3$ and $n = p+1$ or $n = p+2$ \cite{Co}

The paper is organized as follows: in Section 2, we review some known results about bismash product Hopf 
algebras constructed from factorizable groups and the description of the simple modules for such Hopf algebras, 
as well as basic facts about Frobenius-Schur indicators. In Section 3 we prove that the Hopf algebra $H_n$ 
described above is totally orthogonal. In Section 4 we prove our results about more general bismash products, 
and in Section 5 we apply the results of Section 4 to prove our results about $J_n$. We also compute some 
explicit examples. 

Throughout, $H$ will be a finite dimensional semisimple Hopf algebra
over an algebraically closed field $k$ of characteristic $0,$ with
comultiplication $\Delta:H \to H \otimes H$ given by $\Delta(h)=\sum
h_1\otimes h_2$, counit $\epsilon:H\to k$ and antipode $\al:H\to H$.
In particular we know from \cite{LR2} that $H$ is cosemisimple and
from \cite{LR1} that $\al^2=id.$ By Maschke's theorem there exists a
unique integral $\Lambda \in \int _H$ with $\epsilon (\Lambda)=1$.


\section{Extensions arising from factorizable groups and their representations}

The Hopf algebras we consider here were first described by Takeuchi \cite{T}, constructed from 
what he called a matched pair of groups. These Hopf algebras can also be constructed from a 
factorizable group, an old topic in finite groups. We will use both notions here, depending on which 
is more convenient. Throughout, we assume that $F$ and $G$ are finite groups.

\begin{df}
{\rm A \textit{matched pair of groups} $(F,G,\rhd ,\lhd)$ is a pair of
groups $F,G$ with actions $G\overset {\lhd}{\leftarrow} G \times
F\overset{\rhd}{\rightarrow}F$ such that for all $a,b\in F,\ x,y\in G$} 
$$x\rhd ab=(x\rhd a)((x\lhd a)\rhd b),$$ $$xy\lhd a=(x\lhd(y\rhd a))(y\lhd a).$$
\end{df}

\begin{df}
{\rm A group $L$ is called \textit{factorizable} into subgroups $F,G \subset L$ if $FG= L$ and 
$F\cap G = 1$; equivalently, every element of $L$ may be written uniquely as a product $l = ax$ 
with $a\in F$ and $x\in G$}. 
\end{df}

It is easy to see that a factorizable group gives rise to a matched pair, if we define the mutual actions by 
$xa = (x\rhd a)(x\lhd a)$ for $a\in F$, $x\in G$. That is, let $x\rhd a\in F$ and $x\lhd a\in G$ be the (necessarily 
unique) elements of $F$ and $G$ whose product is $xa$. 

Conversely, if  $(F,G,\rhd ,\lhd)$ is a matched pair, then we can define a group $L = F\bowtie G$ with underlying 
set $F\times G$ and multiplication 
$$(a,x)(b,y) = (a (x\rhd b),(x\lhd b)y)$$
for $a,b\in F$ and $x,y\in G$.  This group $F\bowtie G$ is factorizable into the subgroups $F' = 
F\times\{ 1\} \cong F$ 
and $G' = \{ 1\}\times G\cong G$. 

While the actions $\rhd$ and $\lhd$ of $F$ and $G$ on each other are only module actions, they induce 
actions of $F$ and $G$ as automorphisms of the dual algebras $k^G$ and $k^F$. To see this, 
choose a basis $\{p_x\; |\; x\in G\}$ of $k^G$ dual to the elements $x\in G$; similarly let 
$\{p_a\; |\; a\in F\}$ be the basis dual to the elements $a\in F$. Then the new actions are given by 
$$a\cdot p_x := p_{x\lhd a^{-1}}$$ and $$x\cdot p_a := p_{x\rhd a}, $$
for all $a\in F$, $x\in G$,

These actions enable us to construct the {\it bismash product} Hopf algebra $H= k^G\# kF$ associated 
to the matched pair $(F,G,\rhd ,\lhd)$ (see \cite{Ma2}, \cite{Ma3} for more details). As a vector space, 
$H = k^G\otimes kF$, 
with basis $\{p_x\# a \; | \; x\in G, a\in F\}$. The algebra structure is the usual smash product, given by 
\begin{equation}\num \label{mult}
(p_x\# a)(p_y\# b)= p_x (a\cdot p_y)\# ab = p_x p_{y\lhd a^{-1}} \# ab = \delta _{y,{x\lhd a}}p_x\# ab.
\end{equation}
Similarly the dual algebra $H^* = k^F \# kG$ is also a smash product, using the action of $G$ of $k^F$; 
the algebra structure of $H^*$ dualizes to give the coalgebra structure of $H$. More explicitly, this is 
given by 
\begin{equation}\num \label{comult}
\Delta(p_x\# a)=\sum_{y\in G}p_{xy^{-1}}\#( y\rhd a)\otimes p_y\# a.
\end{equation}
Finally $H$ has counit $\epsilon (p_x\# a)=\delta _{x,1}$ and antipode  $\al(p_x\# a)= 
p_{(x\lhd a)^{-1}}\# (x\rhd a)^{-1}$.

The bismash product $H$ is an example of an abelian extension. 
That is, there is an exact sequence of Hopf algebras 
\begin{equation*}
k^G\overset{i}{\hookrightarrow }H\overset{\pi }{\twoheadrightarrow}kF.
\end{equation*}
One could construct a more general extension from $k^G$ and $kF$, a {\it bicrossed product}, in 
which the actions are twisted by cocycles $\sigma: F\times F \to k^G, \tau: G\times G \to k^F$. 
In fact bicrossed products give all extensions:  

\begin{prop} \cite{Ma2} \cite{Ma3}
Let $(F,G,\rhd, \lhd)$ be a matched pair. Every extension of $kF$ by
$k^G$ is equivalent to an extension of the form
\begin{equation*}
k^G\overset{i}{\hookrightarrow } k^G*_{\sigma, \tau}kF \overset{\pi
}{\twoheadrightarrow}kF.
\end{equation*}
\end{prop}

The set of equivalence classes of extensions associated to the matched pair 
is denoted by $Opext(kF,k^G)$. It has a group structure, with identity given 
by the bismash product $k^G\# kF$ as above; it is the extension with trivial 
cocycles \cite{Ma2}, \cite{Ma3}. It turns out that in the examples of interest in this paper, 
we only need look at bismash products:

\begin{ex} \label{Sn} {\rm Let $S_n$ be the symmetric group of degree n, consider 
$S_{n-1} \subset S_n$ by letting any $\sigma\in S_{n-1}$ fix $n$, and let $C_n = \langle z \rangle$, 
the cyclic subgroup of $S_n$ generated by the $n$-cycle 
$z =  (1,2, \ldots,n)$. Then $S_n = S_{n-1}C_n = C_nS_{n-1}$ shows that $S_n$ is factorizable; equivalently $(S_{n-1},C_n)$ and $(C_n,S_{n-1})$ form matched pairs giving the group $L = S_n.$

In this case, it is shown by Masuoka \cite[Corollary 8.2 and Exercise 5.5]{Ma2} that over 
$k =\mathbb C$, the two groups $Opext(kC_n,k^{S_{n-1}})$ and $Opext(kS_{n-1},k^C_n)$ are 
both trivial. Thus the only Hopf algebras one can obtain in this way from the two factorizations 
of $S_n$ are the two bismash products $H_n=k^{C_n}\# kS_{n-1}$ and 
$J_n=k^{S_{n-1}}\#kC_n \cong H_n^*$.}
\end{ex}

We next review the description of the simple modules over a bismash product. 

\begin{prop} \cite[Lemma 2.2 and Theorem 3.3]{KMM} \label{KMM} Let $H = k^G \# kF$ be a bismash product, as above. 
For the action $x\lhd a$ of $F$ on $G$, fix one element $x$ in each $F$-orbit of $G$, and let $F_x$ 
be the stabilizer of $x$ in $F$. Let $V$ be a simple left $F_x$-module and let $\hat{V} = kF\ot_{kF_x} V$ denote the induced $kF$-module. 

$\hat{V}$ becomes an $H$-module in the following way: for any $y\in G$, $a,b\in F$, and $v\in V$, 
$$(p_y\#a)[b\ot v]=\delta_{y\lhd(ab), x} (ab\ot v).$$
Then $\hat{V}$ is a simple $H$-module under this action, and every simple $H$-module arises 
in this way. 
\end{prop}

Given a simple $H$-module $V$ and its character $\chi$, we define the function 
$$\nu (\chi):=\chi ( \sum \Lambda _1\Lambda _2),$$ 
where  $\Lambda$ is the unique integral in $H$ with $\eps(\Lambda) = 1$ from the introduction. 

The next result shows that the function $\nu(\chi)$ agrees with the 
description of the Frobenius-Schur indicator given in the introduction:

\begin{thm} \label{LM} \cite{LM}
Let $H$ be a semisimple Hopf algebra over an algebraically closed
field $k$ of characteristic not 2. Then for $\Lambda$ and $\nu$ as
above, the following properties hold:
\begin{enumerate}
\item
$\nu (\chi)=0,1$ or $-1$, for all simple characters $\chi.$
\item
$\nu (\chi) \neq 0$ if and only if $V_{\chi} \cong V_{\chi} ^*$,
where $V_{\chi}$ is the module associated to the character.
Moreover, $\nu (\chi)=1$ (respectively $-1$) if and only if $V_{\chi
}$ admits a symmetric  (respectively skew-symmetric) nondegenerate
bilinear $H$-invariant form.
\item
Considering $\al\in End(H)$, $Tr(\al)=\sum_{\chi} \nu (\chi) \chi (1).$
\end{enumerate}
\end{thm}

Throughout, the Frobenius-Schur indicator of a simple module is just the indicator of its irreducible
character.

\begin{cor}\label{sum} For any semisimple Hopf algebra $H$ over an algebraically closed
field $k$ of characteristic not 2, $Tr(\al)=\sum_{\chi} \chi (1)$ $\iff$ $H$ is totally orthogonal; 
that is, all Frobenius-Schur indicators are 1. 
\end{cor}

\begin{proof} This is clear from the last part of Theorem \ref{LM}, since the values of $\nu (\chi)$ 
are 0,1, or -1. 
\end{proof}
 
Note that the definition of $\nu$, Theorem \ref{LM}, and Corollary \ref{sum} agree with what is 
known for a finite group $L$, that is for a group algebra $H = kL$ . 

It is clear from the Corollary that computing $Tr(\al)$ is important for computing indicators. 
For a group algebra $kL$, since the antipode is given by $\al(w) = w^{-1}$ for any $w\in L$, 
clearly $Tr(\al)=|\{l \in L:l^2=1\}|.$ We will show that the same statement holds for bismash products. 

First, for any group $L$, we define $i_L:= |\{l \in L:l^2=1\}|.$

\begin{lem} \label{involutions}
Let $L=FG$ be a factorizable group and $\al$ the antipode of the
bismash product $k^G\# kF$ as above. Then $Tr(\al)=|\{l \in L:
l^2=1\}| = i_L.$
\end{lem}

\begin{proof}
Let $x\in G,\ a\in F.$ Then
$$\al(p_x\#a)=p_{(x\lhd a)^{-1}}\# (x\rhd a)^{-1},$$
so an element contributes to the trace if and only if $$(x\lhd
a)^{-1}=x$$ and $$(x\rhd a)^{-1}=a.$$\\
 If $x=(x\lhd a)^{-1}$ and
$a=(x\rhd a)^{-1}$ then $$xa=(x\lhd a)^{-1}(x\rhd a)^{-1}=((x\rhd a)(x\lhd a))^{-1}=(xa)^{-1}.$$ 
Hence $(xa)^2=1.$ 

Conversely, given $l\in L$ with $l^2=1,$ write $l=xa$ with $x\in G,
a\in F$. Then we have $xa=(xa)^{-1}=a^{-1}x^{-1}.$ Hence $x\rhd
a=a^{-1}$, $x\lhd a=x^{-1}$ and $p_{(x\lhd a)^{-1}}\# (x\rhd
a)^{-1}=p_x \# a,$ proving the statement.
\end{proof}

In particular we will need the number $i_L$ when $L=S_n$. In that case, some facts are 
known about $i_n := i_{S_n}$. It is easy to see that 

\begin{equation}\num \label{recursion}
i_n = i_{n-1} + (n-1)i_{n-2}.
\end{equation}

Moreover, it was shown in \cite{CHM} that as $n\to \infty$, $i_n$ is asymptotic to $\sqrt{n!}$.


\section{The Hopf algebra $H_n=k^{C_n}\# kS_{n-1}$}

In this section we prove that $H_n$ has the same property as does $S_n$ in regards 
to the Frobenius-Schur indicators of its irreducible representations: that is, they are all positive. 

We first need some basic facts about our factorization of $S_n$. Throughout this section, we assume 
that $S_n$ is factored into subgroups $F\cong S_{n-1}$ and $G\cong C_n = \langle z \rangle$, for 
$z=(1\ 2...n)$, as described in Example \ref{Sn}. 

\begin{remark} \cite{Ma2} {\rm For any $\sigma \in S_n$, we may write $\sigma =az^r \in FG,$ where
$r=n-\sigma^{-1}(n)$ and $a=\sigma z^{-r}$. The elements $a\in S_{n-1}$ and $z^r\in C_n$ 
are unique with this property. 

If we begin with $b\in S_{n-1}$ and $x\in C_n$ and let $\sigma = xb$, we have 
$xb= az^r$ for unique $a\in S_{n-1}$ and $r\in \{0,1,\ldots,n-1\}$. Thus the actions defining
the matched pair are given by $x\rhd b=a$ and $x\lhd b=z^r.$}
\end{remark}

We first determine the possible orbits and stabilizers for the action of $S_{n-1}$ on $C_n$.

\begin{lem} \label{orbits1} There are exactly two orbits for the action of $S_{n-1}$ on $C_n$, 
namely $\CO_1 = \{1\}$ and $\CO_z = \{z^r: 1\leq r\leq n-1\}$.

Moreover $F_1 = Stab_{S_{n-1}}(1) = S_{n-1}$ and $F_z = Stab_{S_{n-1}}(z) \cong S_{n-2}.$
\end{lem}

\begin{proof} Choose $a\in S_{n-1}$. Then $za=(z\rhd a)(z\lhd a),$ where by definition of the actions 
$z\lhd a=z^r,$ for $ r=n-a^{-1}($n-1$)$. Thus if we let $a_i=($n-1,n-$i)\in S_{n-1},$ then $z\lhd a_i=z^i$ 
for each $i = 1, \ldots, n-1$. Thus $z\lhd S_{n-1} = C_n-\{1\}.$

Clearly $F_1 = F = S_{n-1}$, and one may check that 

$F_z= \{a\in S_{n-1}: za = (z\rhd a)z \} = \{a\in S_{n-1}: a(n-1)=n-1\} \cong S_{n-2}.$
\end{proof}

\begin{df} {\rm Let $c_n$ be the number of simple $S_n$-modules and let $d_{n,i}$, 
for $i=1, \ldots, c_n$ be their dimensions.}
\end{df}

\begin{lem} \label{dims}
The total number of simple $H_n$-modules is $c_{n-1}+c_{n-2},$
with dimensions $d_{n-1, i},\ i=1,\ldots, c_{n-1}$ and $(n-1)d_{n-2, i},\
i=1,\ldots,c_{n-2}.$
\end{lem}

\begin{proof} We know from Proposition \ref{KMM} that if $\hat{V}$ is a simple
$H_n$-module, then it is induced up from a simple module $V$ of
$kF_x$, for some $x\in C_n.$ Moreover as $x$ ranges over a choice of
one element in each $S_{n-1}$-orbit of $C_n$, we get all simple $H_n$-modules. 
By Lemma \ref{orbits1}, we know there are only two orbits, and so we may choose 
$x=1$ and $x = z$ for the elements in each orbit. 

If $x = 1$, $F_x = S_{n-1}$. Then for any simple $S_{n-1}$-module $V$, we get 
a simple $H_n$-module from the induced $H_n$-module $\hat{V}=
kS_{n-1}\otimes _{kS_{n-1}} V\cong V.$ Thus, exactly $c_{n-1}$ simple $H_n$-modules, with dimensions $d_{n-1,i},\ i=1, \ldots, c_{n-1},$ arise in this way. 

When $x=z$, $F_x = S_{n-2}$. In this case we get a simple $H_n$-simple 
$\hat{V}=kS_{n-1}\otimes _{kS_{n-2}} V$ for each simple module $V$ of $S_{n-2}.$ The 
dimension of each such induced module $\hat{V}$ is $[S_{n-1}:S_{n-2}]dim(V)=(n-1)dim(V),$ 
giving $c_{n-2}$ simple modules of $H_n$ with dimensions $(n-1)d_{n-2,i},\ i=1, \ldots, c_{n-2}.$
This completes the proof.
\end{proof}

We need one more fact about $S_n$. 

\begin{lem} \label{trS_n} $Tr(\al_{S_n}) = i_n=
(n-1)\sum _{i=1} ^{c_{n-2}} d_{n-2,i}+ \sum _{i=1} ^{c_{n-1}}d_{n-1,i}.$
\end{lem}

\begin{proof} We use the fact that $S_n$ has all of its Frobenius-Schur 
indicators equal to 1, and thus by Corollary \ref{sum}, $Tr(\al_{S_n})$ is 
the sum of the degrees of the irreducible characters of $S_n$, for all $n$. Thus, 
using the recursion (\ref{recursion}), 
\begin{eqnarray*}
 Tr(\al_{S_n})&=&i_n = (n-1)i_{n-2}+ i_{n-1}\\
&=&(n-1)Tr(\al_{S_{n-2}})+Tr(\al_{S_{n-1}})\\
&=&(n-1)\sum _{i=1} ^{c_{n-2}} d_{n-2,i}+ \sum _{i=1} ^{c_{n-1}}d_{n-1,i}.
\end{eqnarray*}
\end{proof}

We are now able to prove our main result on $H_n$. 

\begin{thm}
Let $H_n=k^{C_n}\# kS_{n-1}$ be the bismash product associated with
the matched pair $(S_{n-1},C_n)$ arising from the factorizable group
$S_n$. Let $V$ be a simple $H_n$-module with character $\chi$. Then
$\nu (\chi)=1.$
\end{thm}

\begin{proof} \label{H_n} By Corollary \ref{sum}, it suffices to prove that 
$Tr(\al_H)=\sum_{\chi} \chi (1),$ where the sum is over all simple characters of
$H.$ From Lemma \ref{dims} we know that the degrees $\chi (1)$ of
the simple characters of $H_n$ are $d_{n-1,i},\ i=1, \ldots, c_{n-1}$ and 
$(n-1)d_{n-2,i}, \ i=1, \ldots, c_{n-2}.$

On the other hand, by Lemma \ref{involutions}, we know that $Tr(\al_{H_n}) = 
Tr(\al_{S_n}).$ We are now done by Lemma \ref{trS_n}, since $Tr(\al_{S_n})$ is 
exactly the sum of the degrees of the irreducible characters of $H_n$.
\end{proof}


\section{Further results about bismash products}

In this section we prove some general results about the actions $\rhd$ and $\lhd$ of $F$ and $G$ 
on each other in a factorizable group $L = FG$.  We are particularly interested in the stabilizers 
$F_x$, for $x\in G,$ in two cases: either when $F_x = \{1\}$ or when $F_x = F$ (that is when 
$x\in G^F$). We will apply most of these results to study indicators for the Hopf 
algebra $J_n=k^{S_{n-1}}\#kC_n$ in Section 5. 

 \begin{lem} \label{sub} $G^F = \{ x\in G : x\lhd a = x, \; {\rm all}\; a\in F\}$  is a subgroup of $G$ and 
$F^G$ is a subgroup of $F$. In particular, $1\in G^F$.
 \end{lem}

\begin{proof} This follows since $x\lhd a = x$ $\iff$ $xa = (x\rhd a)x$. But then for $x, y\in G^F$, 
$xya = x(y\rhd a)y = x\rhd(y\rhd a) xy$ and thus $xy\lhd a = xy$, for all $a\in F$. 

Similarly, $F^G$ is a subgroup of $F$.
\end{proof}

For a given $x\in G$, the following subset of $F$ will be 
important. We define 
$$F_{x, x^{-1}} := \{ a\in F : x\lhd a = x^{-1}\}.$$ 
Of course if $x^2 = 1$, then $F_{x, x^{-1}} = F_x$.

\begin{lem} \label{inverse}
Let $(F,G,\rhd, \lhd)$ be a matched pair of groups, as above. Let $a\in F$ and
$y\in G.$ Then
\begin{enumerate}
\item
$(y^{-1}\lhd a)^{-1}=y\lhd (y^{-1}\rhd a).$
\item
If $y\in \CO_x$, then $y^{-1}\in \CO_{x^{-1}}.$
\item
$a\in F_{y,y^{-1}} \Leftrightarrow a^{-1}\in F_{y^{-1},y}$ and thus 
$|F_{y,y^{-1}}|=|F_{y^{-1},y}|.$ 
\item
If $a\in F_{y^{-1},y}$ then $(y^{-1}a)^2=(y^{-1}\rhd a)a \in F.$
\end{enumerate}
\end{lem}

\begin{proof} (1) Using the relation $xz\lhd a=(x\lhd (z\rhd a))(z\lhd a)$ with 
$x=y,\ z=y^{-1}$ we get $1 = 1\lhd a=yy^{-1}\lhd a=(y\lhd (y^{-1}\rhd a))(y^{-1}\lhd a).$ 
Thus $(y^{-1}\lhd a)^{-1}=y\lhd (y^{-1}\rhd a).$

(2) If $y\in \CO_x$, then $y=x\lhd c$ for some $c$. Then $y^{-1}=(x\lhd
c)^{-1}=x^{-1}\lhd (x\rhd c),$ by part (1). Thus $y^{-1}\in \CO_{x^{-1}}.$

(3) Note that $y\lhd a=y^{-1} \Leftrightarrow y=y^{-1}\lhd a^{-1}$, so the 
first condition holds. Consequently $|F_{y,y^{-1}}|=|F_{y^{-1},y}|.$

(4) $y^{-1}a=(y^{-1}\rhd a)(y^{-1}\lhd a)=(y^{-1}\rhd a)y$ and right
multiplication by $y^{-1}a$ gives $(y^{-1}a)^2 =(y^{-1}\rhd a)a.$
\end{proof}

\begin{cor} \label{orbitcor} Let $x\in G$ and consider $\CO_x$. Then the following are 
equivalent:

(1)  $x^{-1}\in \CO_x,$ or $|F_{x^{-1},x}|\neq 0$;

(2) $y^{-1} \in \CO_x$ for some $y\in \CO_x$;

(3) $y^{-1} \in \CO_x$ for all $y\in \CO_x.$

\noindent Consequently if $|\CO_x|$ is odd, then $x^{-1}\in \CO_x$ $\iff$ for some $y\in \CO_x$, $y^2 = 1$.
\end{cor}

\begin{proof} The equivalence of (1), (2), and (3) follows from part (2) of Lemma \ref{inverse}. 
\end{proof}

\begin{lem} \label{Fxx'} Consider $x\in G$ with $F_x=\{1\}$. Then $|F_{y,y^{-1}}|\leq 1$, and 
$|F_{y,y^{-1}}| = 1$ for all $y\in \CO_x$ $\iff$ $x^{-1}\in \CO_x.$
\end{lem}

\begin{proof} Assume $F_x=\{1\}$ and that $y\in \CO_x$. Let $a,b\in F_{y,y^{-1}}$ and say
$y=x\lhd c.$ Then
$$x\lhd ca=(x\lhd c)\lhd a=y\lhd a=y^{-1}=y\lhd b=(x\lhd c)\lhd
b=x\lhd cb.$$ Thus $ca(cb)^{-1}\in F_x=\{1\}$ and $a=b,$ so
$|F_{y,y^{-1}}|\leq 1.$ 

For $|F_{y,y^{-1}}|\neq 0$ we need $\CO_x=\CO_{x^{-1}},$ or equivalently (by the previous Corollary), 
$\CO_y=\CO_{y^{-1}}.$
\end{proof}

For any bismash product $H=k^G\# kF$, we know that $\La=\ds{\frac{1}{|F|}}\sum _{a\in F} p_1\# a$ is 
the (unique) integral for $H$ such that $\eps(\La) = 1$. To compute the Schur indicators we need 
$\La^{[2]} = m(\Delta (\Lambda)) = \sum \La_1 \La_2$. Now
$$\Delta (\Lambda)= \ds{\frac{1}{|F|}} \sum _{a\in F} \Delta(p_1\# a)= \frac{1}{|F|} \sum _{a\in F} 
\sum _{y\in G} p_y\# (y^{-1}\rhd a)\otimes p_{y^{-1}}\# a.$$
Thus $\Lambda^{[2]}$ is given by
$$
\La^{[2]}= \ds{\frac{1}{|F|}}\sum_{a\in F} \sum_{y\in G} \delta_{y^{-1},y\lhd (y^{-1}\rhd a)} p_y\# (y^{-1} \rhd a)a.
$$

Now   $y^{-1}=(y^{-1}\lhd a)^{-1}$ $\iff$ $y^{-1}=y\lhd (y^{-1} \rhd a)$ by Lemma \ref{inverse}, and so $\La^{[2]}$ 
becomes
\begin{equation}\num \label{L2}
\La^{[2]}= \frac{1}{|F|}\sum_{y\in G} \sum_{a\in F_{y^{-1},y}} p_y\# (y^{-1}\rhd a)a.
\end{equation}

Now consider an irreducible $H$-module $\hat V$ with character $\hat{\chi}$. Recall from 
Theorem \ref{KMM} that for some fixed $x\in G$, $\hat V$ is induced from an irreducible 
$kF_x$-module $V$ with character $\chi$. To compute the values of $\hat \chi$, we use a 
formula which is a simpler version of \cite[Proposition 5.5]{N2}; it is also similar to the formula 
used in \cite[p 898]{KMM}:

\begin{lem}\label{lemma3} The character $\hat{\chi}$ of $\hat V$ may be computed as follows: 
$$\hat{\chi}(p_y\#a) = \sum_{b^{-1}ab \in F_x} \delta_{y\lhd b,x} \chi (b^{-1}ab),$$
for any $y\in G,\ a \in F$, where the sum is over all $b$ in a fixed set $T_x$ of representatives of the 
right cosets of $F_x$ in $F$, and $\chi$ is the character of $V$.
\end{lem}

\begin{proof} Let $W$ be a basis for $V$. Then $\{b\otimes p_x\otimes v:\ b\in T_x,\ v\in W\}$ 
is a basis for $\hat{V}$ (the $p_x$ is inserted to emphasize that the representation depends on $x$). 
Then the action of $p_y\#a$, as in Proposition \ref{KMM}, is given on the elements of this basis by
$$(p_y\#a)[b\otimes p_x\otimes v]=(p_y\# ab)[1 \ot  p_x \ot v] =
(1\#b)(p_{y\lhd b}\# b^{-1}ab)[1\ot p_x\otimes v],$$
using Equation \ref{mult}.
In order to compute the trace of this action, we only need to consider those basis vectors 
$[b\otimes p_x\otimes v]$ for which $b^{-1}ab\in F_x$, since for such $b$ we have 
$$(p_y\#a)[b\otimes p_x\otimes v] = \delta_{y\lhd b,x}b\otimes p_x\otimes (b^{-1}ab)v.$$
This proves the lemma. 
\end{proof} 

Thus for a fixed $x\in G,$ if $\chi$ is the character of $V$ and $\hat{\chi}$ is the character of $\hat{V}$, 
$\hat \chi (p_y\#a)=\chi (b^{-1}ab),$ provided that $b^{-1}ab\in F_x$ and $y\lhd b=x$. 
In particular, in order to have $\hat \chi (p_y\# a) \neq 0$ we need $y\in \CO_x$. This fact together with 
Equation \ref{L2} implies that 
\begin{equation}\num \label{nu}
\hat{\chi}(\Lambda ^{[2]})=\frac{1}{|F|} \sum _{y\in\CO_x} \sum_{a\in F_{y^{-1},y}} \hat{\chi} 
(p_y\#(y^{-1}\rhd a)a).
\end{equation}

\begin{prop}\label{Fx1}
Let $x\in G$ with $F_x=\{1 \}.$ For the unique simple $F_x$-module
$V\simeq k$, the corresponding simple $k^G\# kF$-module $\hat{V}$
has Schur indicator equal to $1$ if and only if $x^{-1}\in
\CO_x$. Otherwise the indicator is $0$.
\end{prop}

\begin{proof} Since $F_x=\{1\},$ by Lemma \ref{Fxx'} $|F_{y^{-1},y}|\leq 1$ for all $y\in \CO_x.$ 
There are two cases, depending on whether or not $x^{-1}\in \CO_x$. 

We first assume $x^{-1}\in \CO_x$. Then for all $y\in \CO_x$, $|F_{y^{-1},y}| = 1$ by 
Lemma \ref{Fxx'}. For each such $y$, let $a\in F_{y^{-1},y}$ then $y\lhd(y^{-1}\rhd a) = 
(y^{-1}\lhd a)^{-1}=y^{-1}$ implies that $y^{-1}\rhd a \in F_{y,y^{-1}}=\{a^{-1} \},$ using Lemma 
\ref{inverse}. Then $(y^{-1}\rhd a)a=1$ and so $\sum _{a\in F_{y^{-1},y}}p_y\# (y^{-1}\rhd a)a = 
p_y\# 1$. Thus Equation \ref{nu} becomes 
$$\nu (\hat{\chi}) =  \frac{1}{|F|}\sum _{y\in \CO_x} \hat{\chi} (p_y\# 1) = 
\frac{1}{|F|}\sum _{y\in \CO_x} \chi(1)=\frac{1}{|F|}|\CO_x|=1,$$ 
where $|\CO_x|=|F|$ since $F_x=\{1\}.$ 

In the second case, when $x^{-1}\notin \CO_x$, then $F_{y,y^{-1}}= \emptyset$ for all $y\in \CO_x$ 
and so $\nu(\hat{\chi})=0.$
\end{proof}

We now turn to the situation in which $x\in G$ with $F_x=F$; that is, $x\in G^F$. This divides into 
several cases 

\begin{prop} \label{FxF} Consider $x\in G$ with $F_x=F$. For any such $x$, any 
simple module $V$ for $kF_x$ is already an $F$-module, and so $\hat V = V$ 
becomes a simple $k^G\#kF$-module, with the additional $k^G$-action. 
There are three possible cases: 

(1) Assume that $x=1$. Then the value of the Frobenius-Schur indicator of $\hat V$ is 
exactly the same as the indicator for $V$ as a simple $kF$-module.

(2) Assume that $x^2\neq 1$. Then the value of the Frobenius-Schur indicator of $\hat V$ 
is $0$.

(3) Assume that $x^2 = 1$ but $x\neq 1$. Suppose in addition that $F=C_{p^r}$, where $p$ 
is an odd prime and $r\in \mathbb{N}$. Then the value of the Frobenius-Schur indicator of 
$\hat V$ is $1$.

\end{prop}

\begin{proof} (1) 
Since $F_{1,1}=F_1=F$, given a simple $kF_1$-module $V$ with
character $\chi$ we have from Equation \ref{nu} that 
$$\nu (\hat{\chi})=\frac{1}{|F|}\sum_{a\in F} \hat{\chi} (p_1\# (1\rhd a)a)=
\frac{1}{|F|}\sum_{a\in F} \chi(a^2)=\nu (\chi).$$

(2) Since $F_x=F$, $x\lhd a=x,\ \forall a\in F.$ Then
$\mathcal{O}_x=\{x\}$ and given that $x\neq x^{-1}$,
$F_{x^{-1},x}=\emptyset.$ Thus for the only element in the orbit of
$x$, namely $x$, there are no elements in $F$ that invert it; 
therefore again by Equation \ref{nu}, $\nu (\hat{\chi})=0.$

(3)
Consider $\phi: F \to F$ given by $\phi (a)=x\rhd a.$ Since $F_x=F,$
$$xa=(x\rhd a)(x\lhd a)=(x\rhd a)x$$ and  $$\phi (a)=x\rhd a=xax^{-1}.$$
Therefore $\phi$ is conjugation by $x$, an automorphism of $F$ of
order two. When $F=C_{p^r},\ Aut(F)$ is cyclic, and it therefore contains 
a unique element of order $2.$ But
$\psi: F\to F,\ \psi (a)=a^{-1}$ is such an automorphism, thus $\phi
= \psi.$ Finally, $$\nu (\hat{\chi})=\frac{1}{|F|} \sum _{a\in F}
\hat{\chi} (p_x\# (x\rhd a)a)=\frac{1}{|F|} \sum _{a\in F} \chi
((xax)a)=\frac{1}{|F|} \sum _{a\in F} \chi (1)=\chi (1).$$
\end{proof}

Note that as long as $G$ contains an element $x$ as in (2), $k^G\#kF$ admits 
simple modules with indicator $0$. This is the case for $J_n$, as we will see in the next 
section. It did not happen for $H_n$, since in that case, only $x =1$ had $F_x=F$, 
and so we were always in case (1) and $F = S_{n-1}$ has all indicators 1. 

Concerning part (3), we do not know what happens when $x^2 = 1$ but $x\neq 1$, if 
no assumptions are made about $F$. 

We are able to obtain one general result.

\begin{thm} \label{Cp} Let $L = FG$ be a factorizable group with $F \cong C_p$, for an odd prime $p$, 
and let $H = k^G \# kF$ be the bismash product. For a given $x\in G$ and its $F$-orbit $\CO_x$, 
let $\hat{V}$ be the simple $H$-module coming from a simple $kF_x$-module $V$. Then:

(1) if $x = 1$, then $\nu(\hat{V}) = 1$ only if $V$ is the trivial representation; otherwise $\nu(\hat{V}) = 0;$

(2) if $x\neq 1$, then $\nu(\hat{V}) = 1$ if and only if there exists some $y\in \CO_x$ such that $y^2=1$; otherwise $\nu(\hat{V}) = 0.$ 
\end{thm} 

\begin{proof} First note that there are only two possibilities for $x$, that is either $F_x = \{1\}$ (in which 
case dim$\hat{V} = p$) or $F_x = F$ (in which case dim$\hat{V} = 1$). 
 
If $F_x = \{1\}$ then $|\CO_x|=p$, an odd number, and so by Corollary \ref{orbitcor}, $x^{-1}\in \CO_x$ 
if and only if there is some $y\in \CO_x$ with $y^2=1$. If some such  $y\in \CO_x$, then the indicator 
of $\hat V$ is 1, by Proposition \ref{Fx1}. 
If no such $y\in \CO_x$ then the indicator is 0.

If $F_x = F$, first assume that $x = 1$, so that $\CO_1=\{1\}$. By \ref{FxF}, 
$\nu(\hat{V}) = \nu(V)$. But these values are all 0, except for the trivial representation, since 
$F$ has odd order. This proves (1) in the theorem. 

If $x\neq 1$, then \ref{FxF} (2) and (3) finish the proof. 
\end{proof} 

In fact we can compute the number of each type of representation more precisely. 

\noindent {\bf I: $\hat{V}$ has dimension 1}. Here $F_x = F$, so that $x\in G^F$. Then either 

(a) $x= 1$, in which case we have the trivial module $\hat{V}_0$ with $\nu(\hat{V}_0) = 1$, plus 
$p-1$ other 1-dim simple modules $\hat{V}$, all with indicator 0; 

(b) $x$ has order 2, in which case we have $p$ simple modules with $\nu(\hat{V}) = 1$;

(c) $x$ has order $> 2$, in which case we have $p$ simple modules with $\nu(\hat{V}) = 0$.

\noindent {\bf II. $\hat{V}$ has dimension $p$}. Here $F_x = \{1\}$, and there are two cases for $x$:

(d) $\CO_x$ contains an element of order 2. Then there is one simple module $\hat{V}$, with $\nu(\hat{V}) = 1;$ 

(e) $\CO_x$ contains no elements of order 2. Then there is one simple module, with $\nu(\hat{V}) = 0.$ 

\vs

Let $m_{L,1}$ be the number of $\CO_x \subset G$ of type (d) and let $m_{L,0}$ be the number of 
$\CO_x \subset G$ of type (e). Also, recall from Section 2 that for any group $L$, $i_L$ is the number 
of $w\in L$ such that $w^2 = 1$. Then we have

\begin{cor} \label{count} Let $L = FG$ with $F=C_p$, and let $H = k^G \# kF$ be as in Theorem \ref{Cp}. Then 

(1) $i_L = 1 + p(i_{G^F} -1) + p m_{L,1}$, and 

(2) $m_{L,0} + m_{L,1} = \ds\frac{|G|-|G^F|}{p}$.
\end{cor} 

\begin{proof} By  Lemma \ref{involutions}  and Theorem \ref{LM},  
$$Tr_H(\al)=i_L = \sum_{\chi} \nu (\chi) \chi (1).$$ 
Since $\chi (1) =$ dim$(\hat{V})$ and $\nu(\chi)$ = 0 or 1 by Theorem \ref{Cp}, we must only consider 
those $\hat{V}$'s of type (b) and (d), together with the trivial module. Thus 
$Tr_H(\al) = 1 + p(i_{G^F} -1) + p m_{L,1}$.

The second formula follows by counting the number of orbits $\CO_x$ of size $p$. This is 
just the number of elements in $G-G^F$ divided by $p$. 
\end{proof}


\section{The dual Hopf algebra $J_n=k^{S_{n-1}}\#kC_n$}

We will now apply the results of Section 4 to study the representations of $J_n=k^{S_{n-1}}\#kC_n.$
To do this, we will view $S_n$ as factored into the subgroups $F = C_n$ and $G = S_{n-1}$; 
that is, we use the matched pair $(C_n,S_{n-1})$, which is reversed from what we did for $H_n$. 
In order to apply Theorem \ref{Cp} and its Corollary, we first need to determine the invariant subgroup 
$G^F = (S_{n-1})^{C_n}$ and its number of elements of order 2. 
We write $C_n = \langle a \rangle$, for $a=(1, 2, \ldots, n)$; $x,y$ will denote elements of $G = S_{n-1}.$ We recall from 
Section 2 that $xa^i = (x\rhd a^i)(x\lhd a^i)$.

\begin{prop} \label{S^C} $(S_{n-1})^{C_n}\cong \BZ_n^*$, the multiplicative group of units in $\BZ_n$. 
Thus $|(S_{n-1})^{C_n}| =\Phi(n).$
\end{prop}

\begin{proof} Choose $x\in G$. Then $x\in G^F$ $\iff$ $x\lhd a^i = x$, for all $i$ $\iff$ $x\lhd a= x$ $\iff$ 
$xa = (x\rhd a) x$ $\iff$ $x\rhd a = xa x^{-1}$. Since $x\rhd a\in F =  \langle a \rangle$, it must be that 
$x\rhd a= a^r$, for some $r \in \{0,1,\ldots,n-1\}$. But also $x a x^{-1}$ is an $n$-cycle in $S_n$. Thus 
$$ a^r = xa x^{-1}= ( x(1), x(2), \ldots,x(n-1),n).$$
Now $a^r$ has order $n$ $\iff$ gcd$(r,n)=1$, so there are at most $\Phi(n)$ such equations  $a^r = xa x^{-1}$. 
Each such equation determines one (and only one) $x = x_r$, since $x$ is determined by the values 
$\{x(1), x(2), \ldots, x(n-1)\}$, and so $x$ is determined by $r$. 

Conversely, given $r$ with gcd$(r,n)=1$, then $a^r$ is also an $n$-cycle. Rearrange $a^r$ (if necessary) so 
that $n$ appears last. Then we may {\it define} $x\in S_{n-1} $ using the first $n-1$ entries of  $a^r$, in order. 
Thus there exists such an $x = x_r$. One can see that $x_r\in G^F$ by reversing the computation at the 
beginning of the proof. 

We may then define a map $f: (S_{n-1})^{C_n} \to \BZ_n^*$, via $x = x_r \mapsto r$; clearly $f$ is 
a group homomorphism, and thus an isomorphism. The Lemma follows. 
\end{proof}

The proof of the proposition shows how to find $(S_{n-1})^{C_n}$ explicitly. For any $r$ relatively 
prime to $n$, 
$$a^r = (r, \overline{2r}, \overline{3r}, \ldots, \overline{(n-1)r}, n),$$
where $\overline{ir}$ is the smallest positive remainder upon dividing $ir$ by $n$. Thus as a permutation in $S_{n-1}$, 
$x_r(i) = \overline{ir}$, for $i = 1, \dots, n-1.$

\begin{cor} \label{cyclic} 
When $p$ is an odd prime, the group $(S_{p-1})^{C_p}$ is cyclic of order $p-1$, generated 
by a $(p-1)$-cycle. 
\end{cor}

\begin{proof} By the Proposition,  $(S_{p-1})^{C_p}\cong \BZ_p^*$, which is cyclic of order $p-1$. To see that it is 
generated by a $p-1$-cycle, we use the description of $x_r$ before the Corollary. It follows that 
$x_r^k(i) \equiv ir^k$ (mod $p$). Choose $r$ to be a primitive root for $p$; then $\{ 1, r, r^2, \ldots, r^{p-2}\}$ are 
all distinct (mod $p$) and so $x_r = (1, r, r^2, \ldots, r^{p-2})$ is a $(p-1)$-cycle.
\end{proof}

We now wish to apply Corollary \ref{count} to $J_p=k^{S_{p-1}}\#kC_p.$ We slightly simplify our notation: let $m_{p,0}=m_{S_p,0}$, $m_{p,1}= m_{S_p,1}$, and $i_p = i_{S_p}$.

\begin{thm}\label{Jp} Let  $J_p=k^{S_{p-1}}\#kC_p.$ 
Then the numbers $m_{p,0}$ and $m_{p,1}$ in Corollary \ref{count} are completely determined, as follows: 

(1) $m_{p,1} = \ds\frac{i_p-1}{p} -1 $, and 

(2) $m_{p,0} + m_{p,1} = (p-1)\ds\frac{(p-2)! -1}{p}$.

Thus $J_p$ has $p(p-1)$ simple modules of dim 1, with exactly $p+1$ of them having indicator +1 and 
the others having indicator 0. 

All other simple $J_p$-modules have dim $p$: there are $(p-1)\ds\frac{(p-2)! -1}{p}$ of these modules altogether, 
of which $m_{p,1}$ have indicator +1 and the rest indicator 0.
\end{thm}

\begin{proof} By Corollary \ref{cyclic}, $(S_{p-1})^{C_p}$ is cyclic of order $p-1$, and so has a unique 
element of order 2. Thus $i_{G^F} = 2$, so by Corollary \ref{count}, $i_p = 1 + p + p m_{p,1}$ . Equation (1) now follows. 

To show (2), 

$m_{p,0} + m_{p,1} = \ds\frac{|G|-|G^F|}{p} = \ds\frac{(p-1)!-(p-1)}{p}= (p-1)\ds\frac{(p-2)! -1}{p}.$

The number and indicators of the simple modules are now described by Theorem \ref{Cp} and Corollary \ref{count}.
\end{proof}

In fact as $p$ gets very large, the representations with indicator 0 (that is, the non-self-dual representations) predominate. 

\begin{cor} \label{lim} Let $M_{p,1}$ be the number of simple modules for $J_p$ with indicator 1, and let $M_p$ be the total number of simple modules for $J_p$.  Then 
$$\ds\frac{M_{p,1}}{M_p}= \ds\frac{i_p + p^2-1}{(p-1)! + (p^2-1)(p-1)}$$
and is asymptotic to  $\ds\frac{1}{\sqrt{(p-2)!}}$. Thus 
$\lim_{p \to \infty} \ds\frac{M_{p,1}}{M_p} = 0.$ 
\end{cor} 

\begin{proof} From Theorem \ref{Jp}, we have that $m_{p,1} =  \ds\frac {i_p-1}{p} -1.$ Thus 

$M_{p,1} = \ds\frac {i_p -1}{p} -1 + (p+1) = \ds\frac{i_p-1+p^2}{p}$. 

The total number of simple modules is 

$$M_p=  p(p-1) + \ds\frac{p-1}{p} [(p-2)! -1] = \ds\frac{(p-1)! + (p^2-1)(p-1)}{p}.$$
Cancelling $p$ from both $M_{p,1}$ and $M_p$ we obtain the desired ratio. 
We now use the result of \cite{CHM}, mentioned at the end of Section 2, that $i_n$ is asymptotic 
to $\sqrt{n!}$. It follows that $\ds\frac{M_{p,1}}{M_p}$ is asymptotic to $ \ds\frac{1}{\sqrt{(p-2)!}}.$ 
Thus $\lim_{p \to \infty} \ds\frac{M_{p,1}}{M_p} =  0.$
\end{proof}

We use Theorem \ref{Jp} to completely describe the values of the indicators of the simple 
$J_p$-modules, in the two cases $p=5$ and $p=7$. 

\begin{ex} \label{J5} {\rm The case $J_5=k^{S_4}\#kC_5.$ 

\vs
We first need $i_5$, the number of $w\in S_5$ such that $w^2 = 1$: one can easily check that 
$i_5 = 26$. It then follows by (1) of Theorem \ref{Jp} that $m_{5,1} = 4$. Also, $(S_{4})^{C5}$ has 
order 4, generated by a 4-cycle, and $(p-1)\ds\frac{(p-2)! -1}{p} = 4 \ds\frac{(3)! -1}{5} = 4$. 
Thus by (2) of Theorem \ref{Jp}, $m_{5,0} = 0$. 
Thus $J_5$ has 20 simple modules of dim 1, of which 6 have indicator +1, and $J_5$ has 4 simple modules of dim 5, 
all of which have indicator +1. 

The ratio in Corollary \ref{lim} is $\ds\frac{M_{5,1}}{M_5} = \ds\frac{6 + 4}{24}=\ds\frac{10}{24}\sim .417$}
\end{ex}

In the table below, we give a complete list of the different orbits, given by the action of $C_5$ on $S_4$. 
The table was computed using Maple, in particular some of the techniques in \cite{LMS}.

\begin{table}[h]
\caption{The action of $C_5$ on $S_4$} \centering
\begin{tabular}{c|c c c c c|c}
\hline $\lhd$& 1& (12345)& (13524)& (14253)& (15432)& $F_x$\\
[0.5ex] \hline

x=1& 1& 1& 1 &1 &1& $C_5$\\
x=(14)(23)& (14)(23)& (14)(23)& (14)(23)& (14)(23)& (14)(23)& $C_5$\\
x=(1243)& (1243)& (1243)& (1243)& (1243)& (1243)& $C_5$\\
x=(1342)& (1342)& (1342)& (1342)& (1342)& (1342)& $C_5$\\ 
x=(12)& (12)& (1432)& (1234)& (34)& (23)& $\{1\}$\\
x=(13)& (13)& (1423)& (14)& (1324)& (24)& $\{1\}$\\
x=(123)& (123)& (243)& (132)& (13)(24)& (234)& $\{1\}$\\
x=(124)& (124)& (12)(34)& (143)& (134)& (142)& $\{1\}$\\[1ex]

\hline
\end{tabular}
\label{5table}
\end{table}

We observe directly that $(S_4)^{C_5} =  \{ e, (1243), (14)(23), (1342)\} = \langle (1243) \rangle$, 
with a unique element of order 2. $x_2 =(1243)$ is precisely the generator obtained in the Corollary, 
using that $r=2$ is a primitive root of 5. Also $x=(14)(23)$ is the unique element of order 2 in $C_5$, 
and so Theorem \ref{Cp} gives $\nu (\hat{\chi})=\chi (1)$ for $\hat V$ coming from $\CO_x$. 
But $F_x= C_5$, so all irreducible representations of $F_x$ are $1$-dimensional. Thus the five 
irreducible representations corresponding to the stabilizer of $(14)(23)$ all have indicators $+1$.

For $x=(1243)$ and $x=(1342)$, $F_x=C_5$ and $x^2\neq 1$ 
so $\nu(\hat{\chi})=0$, for all characters $\hat{\chi}$ corresponding to their orbits. For $x=1$, we 
know that the Schur indicators are those of $F_1=C_5$. Since $5$ is odd, the only self dual 
irreducible representation of $C_5$ is the trivial one. Thus the five linear irreducible
representations of $C_5$ have indicators $1,0,0,0,0.$

For the $5$-dimensional representations, we can see from the table that for the four elements 
$x=(12),\ (13),\ (123),\ (124)$ with $F_x=\{1\},$ each of their orbits $\CO_x$ contains an 
element of order 2. Thus $\nu(\hat{\chi})=1$ for all the corresponding characters.

\begin{ex} \label{J7} {\rm The case $J_7=k^{S_6}\#kC_7.$ 

\vs

Again we need $i_7$: one may check that $i_7 = 232$. Thus (1) of Theorem \ref{Jp} gives $m_{7,1} = 32$, 
the number of seven-dimensional simple modules with indicator +1. Now 
$m_{7,0} + m_{7,1} = 6\ds\frac{(5)! -1}{7} = 6 \cdot 17 = 102$. Thus there are 70 seven-dimensional simple 
modules with indicator 0. 

Also from Theorem \ref{Jp}, there are $p(p-1)= 42$ simple modules of dim 1, with exactly $8$ of them having 
indicator +1 and the others having indicator 0. 

By the Corollary, $(S_6)^{C_7}$ is cyclic, generated by a 6-cycle.To find a generator, note that $r=3$ 
is a primitive root for 7. Then $x_7 = (1,3,2,6,4,5)\in S_6$ is a generator of $(S_6)^{C_7}$.

The ratio in Corollary \ref{lim} is $\ds\frac{M_{7,1}}{M_7} = \ds\frac{32 + 8}{102+ 42}=\ds\frac{40}{144}\sim .278$

Using Maple, we have also computed all of the orbits of the action of $C_7$ on $S_6$.}
\end{ex}

\begin{ex} \label{J11} {\rm The case $J_{11}=k^{S_{10}}\#kC_{11}.$ 

\vs
In this case we only compute the ratio in Corollary \ref{lim}. Using the recursion \ref{recursion}, 
we see that $i_8 = 764$, $i_9 = 2620$, $i_{10} = 9496$, and $i_{11} = 35,696$. 
Thus the numerator in the ratio is $35,696 + 121 -1= 35, 816$. The denominator is 
$10! + 120(10) = 3628800 + 1200 = 3,630,000.$ Thus} 
$$\ds\frac{M_{11, 1}}{M_{11}} = \ds\frac{35,810}{3,630,000} \sim .00986 \sim .01.$$

\end{ex}

To illustrate what may happen when $n$ is not prime, we look at the case of $n=6$.

\begin{ex} \label{J6} {\rm The case $J_6=k^{S_5}\#kC_6.$} 
\end{ex} 

First, we know from Proposition \ref{S^C} that $S_5^{C_6} \cong (\BZ_6)^* \cong \BZ_2$. 
A direct computation shows that $S_5^{C_6} = \{ e, (15)(24)\}$, using the same method as for 
$ n= 5, 7$.  Using Maple, we computed the action of $C_6$ on some elements of $S_5$, 
to find an orbit with stabilizer $\{1\} < F_x< C_6.$ Since $F_x \neq \{1\}$ or $C_6,$ the methods 
of Section 4 do not apply, making the calculations much more complicated.

\begin{table}[h]
\caption{The action of $C_6$ on part of $S_5$} \centering
\begin{tabular}{c|c c c c c c|c}
\hline $\lhd$& 1& (123456)& (135)(246)& (14)(25)(36)& (153)(264)&
(165432)& $F_x$\\ [0.5ex] \hline

(12)& (12)& (15432)& (12345)& (45)& (34)& (23)& $\{1\}$\\
(13)& (13)& (153)(24)& (15)& (135)(24)& (35)& (24)& $\{1\}$\\
(14)& (14)& (14)(25)& (25)& (14)& (14)(25)& (25)& $<(14)(25)(36)>$\\
[1ex] \hline
\end{tabular}
\label{6table}
\end{table}

When $x=(14),$ table \ref{6table} shows that $F_x=\{1,(14)(25)(36)\},$ a proper subgroup of $C_6$. 
Thus a simple module $V$ for $F_x$ has dim 1, and the corresponding induced module $\hat{V}$ for 
$J_6$ has dimension 3.



\begin{thebibliography}{FGSV}


\bibitem[CHM]{CHM} Chowla, I. N. Herstein, and Moore, On recursions connected with 
symmetric groups, {\it Canadian J. Math} 3 (1951), 328-334.

\bibitem[Co]{Co} Collins, M. J, Some bismash products that are not group algebras, 
{\it J. Algebra} 316 (2007), 297-302.

\bibitem[FGSV]{FGSV}  J. Fuchs, A. Ch. Ganchev, K. Szlach\'{a}nyi and P.
Vecserny\'{e}s, $S_{4}$ symmetry of $6j$ symbols and Frobenius-Schur
indicators in rigid monoidal $C^{\ast }$ categories. {\it J. Math.
Phys.} 40 (1999), no. 1, 408--426.

\bibitem[GW]{GW} L. C. Grove and K. S. Wang, Realizability of representations of finite
groups, {\it J. Pure Appl. Algebra} 54 (1988) 299--310.

\bibitem[GM]{GM} R. Guralnick and S. Montgomery, Frobenius-Schur Indicators for subgroups 
and the Drinfel'd double of Weyl groups, preprint. 

\bibitem[H]{H} J. Humphries, {\it Reflection groups and Coxeter groups},  Cambridge Studies in
Advanced Mathematics 29, Cambridge University Press, 1990.

\bibitem[JL]{JL}  G. James and M. Liebeck, {\it Representations and characters
of groups}, Cambridge Mathematical Textbooks, Cambridge University
Press, Cambridge, 1993. 

\bibitem[KMM]{KMM} Y. Kashina, G. Mason and S. Montgomery, Computing the Frobenius-Schur
indicator for abelian extensions of Hopf algebras, {\it J. Algebra}
251 (2002), 888--913.

\bibitem[KSZ1]{KSZ1} Y. Kashina, Y. Sommerh\"auser, and Y. Zhu, Self-dual modules
of semisimple Hopf algebras,  {\it J. Algebra} 257 (2002), 88--96.

\bibitem[KSZ2]{KSZ2} Y. Kashina, Y. Sommerh\"auser, and Y. Zhu, On higher Frobenius-Schur 
indicators, {\it AMS Memoirs} 181 (2006), no. 855.

\bibitem[LR1]{LR1} R. Larson and D. Radford, Semisimple cosemisimple
Hopf algebras, {\it Am. J. Math.} 109 (1987), 187--195.

\bibitem[LR2]{LR2} R. Larson and D. Radford, Finite-dimensional cosemisimple
Hopf algebras in characteristic 0 are semisimple, {\it J. Algebra} 117 (1988), 267--289.

\bibitem[LMS]{LMS} R. Landers, S. Montgomery and P. Schauenburg,
Hopf powers and orders of some bismash products, {\it J. Pure and
Applied Algebra} 205 (2006), 156--188.

\bibitem[LM]{LM}  V. Linchenko and S. Montgomery, A Frobenius-Schur theorem
for Hopf algebras, {\it Algebras and Representation Theory}, 3
(2000), 347--355.

\bibitem[MaN]{MaN} G. Mason and S-H. Ng, Central invariants and Frobenius-Schur indicators 
for semisimple quasi-Hopf algebras, {\it Advances in Math} 190 (2005), 161-195. 

\bibitem[Ma1]{Ma1}  A. Masuoka, Calculations of some groups of Hopf algebra extensions,
{\it J. Algebra} 191 (1997), 568--588.

\bibitem[Ma2]{Ma2}  A. Masuoka, Extensions of Hopf algebras,
{\it Trabajos de Matem\'{a}tica} 41/99, FaMAF, Universdad Nacional
de C\'{o}rdoba, Argentina, 1999.

\bibitem[Ma3]{Ma3}  A. Masuoka, Hopf algebra extensions and cohomology,
{\it New directions in Hopf algebras}, Mathematical Sciences
Research Institutes Publications, Vol. 3, Cambridge University
Press, Cambridge, 2002, pp. 167--209.

\bibitem[Mo]{Mo}  S. Montgomery, {\it Hopf Algebras and their Actions on Rings},
CBMS Lectures, Vol. 82, AMS, Providence, RI, 1993.

\bibitem[N1]{N1}  S. Natale, On group-theoretical Hopf algebras and exact factorizations of 
finite groups,  {\it J. Algebra} 270 (2003), 199-211.

\bibitem[N2]{N2}  S. Natale, Frobenius-Schur indicators for a class of fusion categories, 
{\it Pacific J. Math} 221 (2005), 363--377.

\bibitem[T]{T}  M. Takeuchi, Matched pairs of groups and bismash products of Hopf algebras, 
{\it Comm. in Algebra} 9 (1981), 84--882.



\end{thebibliography}
\end{document}